\documentclass{elsart}
\usepackage{amssymb}
\newtheorem{theorem}{Theorem}[section]
\newtheorem{lem.}[theorem]{Lemma}
\newtheorem{cor.}[theorem]{Corollary}
\newtheorem{prop.}[theorem]{Proposition}
\newtheorem{definition}[theorem]{Definition}
\newenvironment{proof}{\noindent \bf Proof: \rm}{$ \hspace{\stretch{1}}$ \qed

\vspace{5mm}}
\newcommand{\Hil}[0]{
\mathcal{H} 
}
\newcommand{\MM}[0]{
{\bf M} 
}
\newcommand{\norm}[2]{
\left\| #2 \right\|_{#1}
}
\newcommand{\HS}[0]{
{\mathcal HS}
}
\newcommand{\NN}[0]{
{\mathbb{N}}
}
\newcommand{\BL}[0]{
{\mathcal B}
}
\begin{document}

\sloppy

\begin{frontmatter}
\title{Basic Definition and Properties of Bessel Multipliers}
\author{Peter Balazs}
\address{Austrian Academy of Sciences, Acoustic Research Institute,
         Reichsratsstrasse 17,A-1010 Vienna, Austria}
\date{\today}
\thanks{Manuscript received January 11, 2006; revised February 3, 2005.
        This work was partly supported by the European Union's Human Potential Programe, under contract HPRN-CT-2002-00285 (HASSIP).}
\thanks{Peter.Balazs@oeaw.ac.at}
\begin{abstract}
This paper introduces the concept of Bessel multipliers. These operators are defined by a fixed multiplication pattern, which is inserted between the analysis and synthesis operators.
The proposed concept unifies the approach used for Gabor 
multipliers for arbitrary analysis/synthesis systems, which form Bessel 
sequences, like wavelet or irregular Gabor frames. 
The basic properties of this class of operators are investigated. 
In particular the implications of summability properties of the symbol for the membership of the corresponding operators in certain operator classes are specified. As a special case the multipliers for Riesz bases are examined and it is shown that multipliers in this case can be easily composed and inverted. 
Finally the continuous dependence of a Bessel multiplier on the parameters (i.e. the involved sequences and the symbol in use) is verified, using a special measure of similarity of sequences.
\end{abstract}
\end{frontmatter}
\maketitle
\section{Introduction}
The application of signal processing algorithms like adaptive or time variant filters are numerous \cite{widste1}%
.
If the STFT, the \em Short Time Fourier Transformation \em \cite{Groech1} is used in its sampled version, the Gabor transform, one possibility to construct a time variant filter is the usage of \em Gabor multipliers\em . These operators are a current topic of research \cite{feinow1,doerf1}. For them the Gabor transform is used to calculate time frequency coefficients, which are multiplied with a fixed time-frequency mask and then the result is synthesized. These operators have been already used for quite some time implicitly in engineering applications and recently have been used in signal-processing applications as time-variant filters called \em Gabor filters \em \cite{hlawatgabfilt1}. Recent applications can be found for example in the field of system identification \cite{majxxl1}.

If another way of calculating these coefficients is chosen or if another synthesis is used, many modifications can still be implemented as multipliers. For example it seems quite natural to define 
\em wavelet multipliers\em . Also as irregular Gabor frames get more and more attention \cite{liuwang1}, Gabor multipliers on irregular sets can be investigated \cite{xxlphd1}. As the sampling set, in this irregular case, does not form a lattice, there is no group structure to work with. Therefore it is quite natural for this case to generalize even more and look at multipliers with frames without any further structure.

All these special types of sequences are used in a lot of applications. They have the big advantage, that it is possible to interpret the analysis coefficients. This would also make the formulation of a concept of a multiplier for other analysis / synthesis systems very profitable, like e.g. gammatone filter banks \cite{hartm1}, which are mainly used for analysis based on the auditory system. In \cite{piche1} a gammatone filter bank was used for analysis and synthesis, for the sound separation part a neuronal network creates a mask for these coefficients. This complies with the definition of a multiplier presented here. 

Therefore for Bessel sequences the investigation of operators 
\mbox{$ {\bf M} =  \sum \limits_k m_k \left< f, \psi_k \right> \phi_k $}, 
where the analysis coefficients, $\left< f, \psi_k \right>$, are multiplied by a fixed \em symbol \em $(m_k)$ before resynthesis (with $\phi_k$), is very natural and useful. These are the \em Bessel multipliers \em investigated in this paper. As stated above there are numerous applications of this kind of operators. It is the goal of this paper to set the mathematical basis to unify the approach to them for all possible analysis / synthesis sequences, that form a Bessel sequence.

\section{Main Results}

We will introduce the concept of Bessel multipliers %as a generalization of Gabor multipliers 
and 
will study their basic properties for the first time in an article. 
An important result is dealing with the connection of the symbol, the fixed multiplication pattern, to the operator. Most notably if the symbol is in the sequence spaces $l^\infty, c_0, l^2$ or $l^1$ respectively, then the multiplier is a bounded, compact, trace class or Hilbert-Schmidt operator respectively. We will also prove
that for Riesz bases the Bessel multipliers behave 'nicely', most importantly that the mapping of the symbol to the operator is an injective one. The last result states that the Bessel multiplier depends continuously on the symbol and on the
involved Bessel sequences (in a special sense). For this result the investigation of the perturbation of Bessel sequences is important. This topic is given some thought right after the introduction.

This article is organized as follows: Section \ref{sec:notandrefs1} will first fix some notations and review basic facts in some detail.
In Section \ref{sec:perturbfram0} we are going to present results on the perturbation of Bessel sequences, frames and Riesz bases needed in Section \ref{sec:framchangincred0}. 
Section \ref{sec:frammult0} will give the basic definition and preliminary results for Bessel multipliers.
In Section \ref{sec:frammulprop} we will look at the influence the symbol has on the operator and investigate further properties.
Section \ref{sec:rieszmult1} deals with multipliers for Riesz bases and shows that in this case they behave 'nicely' in many ways.
In Section \ref{sec:framchangincred0} the influence of "small" changes of the parameters on the operator is examined. The paper is finished with Section \ref{sec:conclpersp0}, Perspectives.

This article is based on parts of %the author's PhD thesis 
\cite{xxlphd1}. Some straightforward and easy proofs can be found there and are not given here.

\section{Notation and Preliminaries} \label{sec:notandrefs1}

In this Section basic notation and preliminary result are collected. 
Let $\Hil$ denote a separable Hilbert space. The inner product will be denoted by $\left< ., . \right>$ and will be linear in the first coordinate. Let $\BL(\Hil_1,\Hil_2)$ be the set of all bounded operators from $\Hil_1$ to $\Hil_2$. With the \em operator norm \em, 
$\norm{Op}{O} = \sup \limits_{\norm{\Hil_1}{x} \le 1} \left\{ \norm{\Hil_2}{O (x)} \right\}$,
 this set forms a Banach algebra. By $O^*$ we denote the \em adjoint operator\em. 
 For more details on Hilbert space respectively operator theory see \cite{conw1}.

Recall that an operator $T \in \BL(\Hil_1,\Hil_2)$ is called \em compact\em , if $\overline{T(B_1)}$ is compact with $B_1$ being the unit ball. 
We know that $T$ is compact, if and only if there exist a sequence $T_n \in \BL(\Hil_1,\Hil_2)$ with finite rank, such that $\norm{Op}{T_n - T} \rightarrow 0$ for $n \rightarrow \infty$.
Special classes of compact operators we are using are the trace class (respectively Hilbert-Schmidt class ($\HS$)) operators, which are operators, where the singular values are summable (respectively square summable) with the respective norms $\norm{trace}{.}$ and $\norm{\HS}{.}$. For details see \cite{schatt1}, \cite{wern1} or \cite{xxlphd1}.
We will be using the following special operator:
\begin{definition} \label{sec:kronop1} 
Let $f \in \Hil_1$, $g \in \Hil_2$ then define the \em (inner) tensor product \em as the operator from $\Hil_2$ to $\Hil_1$ by
$ %
\left( f \otimes_{i} \overline{g} \right) (h) = \left< h, g \right> f $ 
\end{definition} 
For this operator we know  
\cite{schatt1} 
that it is a bounded linear operator from $\Hil_2$ to $\Hil_1$ with $\norm{Op}{f \otimes_i \overline{g}} = \norm{\Hil_1}{f} \cdot \norm{\Hil_2}{g}$. The last equality is also true for $\norm{trace}{.}$ and $\norm{\HS}{.}$.

\subsection{Frames and Bases} 
For more details on this topic  see e.g. \cite{ole1} or \cite{Casaz1}. 
\begin{definition}\label{sec:framprop1} 
A sequence $(\psi_k)$ is called a 
\em frame \em for the Hilbert space $\Hil$, if constants $A,B > 0$ exist, such that 
$A \cdot \norm{\Hil}{f}^2 \le \sum \limits_k \left| \left< f, \psi_k \right> \right|^2 \le B \cdot  \norm{\Hil}{f}^2  \ \forall \ f \in \Hil$. \\
$A$ is a \em lower\em , $B$ an \em upper frame bound\em .
If the bounds can be chosen such that $A=B$ the frame is called \em tight\em . A sequence is called \em Bessel sequence \em if the right inequality above is fulfilled.
\end{definition}
The index set will be omitted in the following, if no distinction is necessary. 
The optimal bounds $A_{opt}, B_{opt}$ are the biggest $A$ and smallest $B$ that fulfill the corresponding inequality.
\begin{lem.} \label{sec:frambound2} Let $(\psi_k)$ be a Bessel sequence for $\Hil$. Then 
$ \norm{\Hil}{\psi_k} \le \sqrt{B}$.
\end{lem.}

\begin{definition} For a Bessel sequence $(\psi_k)$ let 
\begin{itemize}
\item $C_{( \psi_k )} : \Hil \rightarrow l^2 ( K )$ be the \em analysis  operator \em
$C_{( \psi_k )} ( f ) = \left( \left< f , \psi_k \right> \right)_k$, 
\item $D_{( \psi_k )} : l^2( K ) \rightarrow \Hil $ be the  \em synthesis operator \em
$D_{( \psi_k )} \left( \left( c_k \right) \right) = \sum \limits_k c_k \cdot \psi_k$ and
\item $S_{( \psi_k )} : \Hil  \rightarrow \Hil $ be the \em (associated) frame  operator \em
$S_{( \psi_k )} ( f  ) = \sum \limits_k  \left< f , \psi_k \right> \cdot \psi_k .$
\end{itemize}
\end{definition}
If there is no chance of confusion, we will omit the index, so e.g. write $C$ instead of $C_{( \psi_k )}$. 
These operators have the following properties:
\begin{prop.} 
\begin{enumerate}
\item Let $(\psi_k)$ be a Bessel sequence. Then $C$ and $D$ are adjoint to each other, $D = C^*$ and so $\norm{op}{D} = \norm{op}{C} \le \sqrt{B}$.  The series $\sum \limits_k c_k \cdot \psi_k$ converges unconditionally. 
\item Let $(\psi_k)$ be a frame. $C$ is a bounded, injective operator with closed range. 
\item Let $(\psi_k)$ be a frame. $S = C^*C = DD^*$ is a positive invertible operator satisfying $A I_\Hil \le S \le B I_\Hil$ and $B^{-1} I_\Hil \le S^{-1} \le A^{-1} I_\Hil$.
\end{enumerate}
\end{prop.}

If we have a frame in $\Hil$, we can find an expansion of every member of $\Hil$ with this frame:
\begin{theorem} 
Let $\left( \psi_k \right)$ be a frame for $\Hil$ with  bounds $A$, $B > 0$. Then $ \left( \tilde{g}_k \right) = \left( S^{-1} \psi_k \right)$ is a frame with bounds $B^{-1}$, $A^{-1} > 0$, the so called \em canonical dual frame\em . %\index{frames!canonical dual}.
 Every $f \in \Hil$ has 
 expansions
$ f = \sum \limits_{k \in K} \left< f, \tilde \psi_k \right> \psi_k $
and 
$ f = \sum \limits_{k \in K} \left< f, \psi_k \right> \tilde \psi_k $
where both sums converge unconditionally in $\Hil$.
\end{theorem}

\begin{definition} 
For two sequences $\{ \psi_k \}$ and $\{ \phi_k\}$  in $\Hil$ we call $G_{\psi_k , \phi_k}$ given by $\left( G_{\psi_k, \phi_k} \right)_{jm} = \left< \phi_m , \psi_j \right>$, $j,m \in K$ the \em cross-Gram matrix \em .
If $(\psi_k) = (\phi_k)$ we call this matrix the \em Gram matrix \em $G_{\psi_k}$. 
\end{definition}

We can look at the operator induced by the Gram matrix, defined for $c \in l^2$ formally as
$ ( G_{\psi_k , \phi_k} c )_j = \sum \limits_k c_k \left< \phi_k , \psi_j \right> $. 
For two Bessel sequences it is well defined and bounded with
$\norm{Op}{G_{\psi_k , \phi_k}} \le \norm{Op}{C_{(\psi_k)}} \norm{Op}{D_{(\phi_k)}} \le  B$ .
\begin{definition} Let $(\psi_k)$ be a complete sequence. If there exist constants $A$, $B >0$ such that the inequalities
$ A \norm{2}{c}^2 \le \norm{\Hil}{\sum \limits_{k \in K} c_k \psi_k}^2 \le B \norm{2}{c}^2 $
hold for all $c \in l^2$ 
the sequence $(\psi_k)$ is called a \em Riesz basis\em . 
A sequence $(\psi_k)$ that is a Riesz basis only for $\overline{span}( \psi_k)$ is called a \em Riesz sequence\em . 
\end{definition}
Every subfamily of a Riesz basis is a Riesz sequence. By the lower and upper bounds for a Riesz sequence we will denote the Riesz bounds on the closed span of the elements.
It is evident that Riesz bases are frames and the Riesz bounds coincide with the frame bounds. See Christensen \cite{ole1}.
\begin{theorem} 
Let $( \psi_k )$ be a frame for $\Hil$. Then the following conditions are equivalent:
\begin{enumerate} \label{sec:riesz2}
\item $( \psi_k )$ is a Riesz basis for $\Hil$.
\item The coefficients $( c_k ) \in l^2$ for the series expansion with $( \psi_k )$ are unique, i.e. the synthesis operator $D$ is injective.
\item The analysis operator $C$ is surjective.
\item  There exists sequence, which is biorthogonal to $( \psi_k )$.
\item $( \psi_k )$ and $( \tilde{g}_k )$ are biorthogonal.
\item $( \psi_k )$ is a basis.
\end{enumerate}
\end{theorem}

We also need the following estimation of the norm of the elements:
\begin{cor.} \label{sec:boundofrieszelem1} Let $(\psi_k)$ be a Riesz basis with bounds $A$ and $B$. Then 
$$ \sqrt{A} \le \norm{\Hil}{\psi_k} \le \sqrt{B} \quad \forall \ k \in K $$
\end{cor.}

\section{Perturbation of Bessel sequences} \label{sec:perturbfram0}

For the perturbation results of Bessel multipliers, Section \ref{sec:framchangincred0}, we need some special results on the perturbation of  Bessel sequences. 
The standard question of perturbation theory is whether the Bessel, frame or Riesz properties of a sequence are shared with 'similar' sequences. A well-known result is the following \cite{ole1}: Let $(\phi_k)_{k=1}^\infty$ be a frame for $\Hil$ with bounds $A$,$B$. Let $(\psi_k)_{k=1}^\infty$ be a sequence in $\Hil$. If there exist $\lambda, \mu \ge 0$ such that $\lambda + \frac{\mu}{\sqrt{A}} < 1$ and
$$ \norm{\Hil}{\sum \limits_k c_k \left(\phi_k - \psi_k\right)} \le \lambda \norm{\Hil}{\sum \limits_k c_k \phi_k} + \mu \norm{l^2}{c}$$
for all finite scalar sequences $c$ (we denote $c \in c_{00}$), then $(\psi_k)$ is a frame with lower bound 
$ A \left( 1 - \left( \lambda + \frac{\mu}{\sqrt{A}} \right) \right)^2$ and upper bound $B \left( 1 + \lambda + \frac{\mu}{\sqrt{B}} \right)^2 $.
Moreover if $(\phi_k)$ is a Riesz bases or Riesz sequence, $(\psi_k)$ is, too. 
This result can be easily formulated for Bessel sequences using only these parts of the proofs in \cite{ole1} Theorem 15.1.1. which apply for Bessel sequences:
\begin{cor.} \label{sec:perturbBessel1} Let $(\phi_k)_{k=1}^\infty$ be a Bessel sequence for $\Hil$. Let $(\psi_k)_{k=1}^\infty$ be a sequence in $\Hil$. If there exist $\lambda, \mu \ge 0$ such that 
$$ \norm{\Hil}{\sum \limits_k c_k \left(\phi_k - \psi_k\right)} \le \lambda \norm{\Hil}{\sum \limits_k c_k \phi_k} + \mu \norm{l^2}{c}$$
for all $(c_k) \in c_{00}$, then $(\psi_k)$ is a Bessel sequence with bound
$ B \cdot \left( 1 + \lambda + \frac{\mu}{\sqrt{B}} \right)^2 $.
\end{cor.}
We can specialize and rephrase these results as needed in the later sections.
 For that let us denote the normed vector space of finite sequences in $l^2$ by $c_{00}^2 = (c_{00}, \norm{2}{.})$. \index{symbols!$c_{00}^2$}
\begin{prop.} \label{sec:perturbframsynth1} Let $(\phi_k)$ be a Bessel sequence, frame, Riesz sequence or Riesz basis for $\Hil$. Let $(\psi_k)$ be a sequence in $\Hil$. If there exists $\mu$ such that 
$$ \norm{c_{00}^2 \rightarrow \Hil}{D_{(\phi_k)} - D_{(\psi_k)}} \le \mu < \sqrt{A} $$
then $(\psi_k)$ shares this property with upper bound
$ B \left( 1 + \frac{\mu}{\sqrt{B}} \right)^2 $
and, if applicable, lower bound
$ A \left( 1 - \frac{\mu}{\sqrt{A}} \right)^2$
and 
$ \norm{l^2 \rightarrow \Hil}{D_{(\phi_k)} - D_{(\psi_k)}} \le \mu  $. 
\end{prop.}
\begin{proof} For every $c \in c_{00}$ we get
$$ \norm{\Hil}{\left( D_{(\phi_k)} - D_{(\psi_k)} \right) c} \le \norm{Op}{D_{(\phi_k)} - D_{(\psi_k)}} \norm{2}{c}  \le \mu \norm{2}{c} .$$
This is just the condition in the perturbation result mentioned above with $\lambda = 0$ and $\mu < \sqrt{A}$, so that $\lambda + \frac{\mu}{\sqrt{A}} < 1$. %\frac{\alpha}{\sqrt{A}}$
Because $(\psi_k)$ is a Bessel sequence, we know that $D_{(\psi_k)} : l^2 \rightarrow \Hil$ is well defined. Because $c_{00}^2$ is dense in $l^2$ we get
$ \norm{l^2 \rightarrow \Hil}{D_{(\phi_k)} - D_{(\psi_k)}} = \norm{c_{00}^2 \rightarrow \Hil}{D_{(\phi_k)} - D_{(\psi_k)}} \le \mu $.
\end{proof}

\begin{cor.} \label{sec:upperboundconverg1} Let $(\phi_k)$ be a Bessel sequence, frame, Riesz sequence respectively Riesz basis and $(\psi_k^{(n)})$ sequences with 
$\norm{c_{00}^2 \rightarrow \Hil}{D_{(\psi_k^{(n)})}-  D_{(\phi_k)}} \rightarrow 0$ 
for $n \rightarrow \infty$. Then there exists an $N$ such that $(\psi_k^{(n)})$ are Bessel sequences, frames, Riesz sequences respectively Riesz bases for all $n \ge N$. For the optimal upper frame bounds we get $B_{opt}^{(n)} \rightarrow B_{opt}$. And
$\norm{l^2 \rightarrow \Hil}{D_{(\psi_k^{(n)})}-  D_{(\phi_k)}} \rightarrow 0$
for $n \rightarrow \infty$.
\end{cor.}
\begin{proof} The first property is a direct consequence from Proposition \ref{sec:perturbframsynth1}.

To show the second part we note that for all $\varepsilon > 0$ there is an $N$ such that for all $n \ge N$, \ 
$\norm{Op}{D_{(\psi^{(n)}_k)}} \le \norm{Op}{D_{(\phi_k)}} + \norm{Op}{D_{(\phi_k)} - D_{(\psi_k^{(n)})}} \le B + \varepsilon $. 
\end{proof}

In Section \ref{sec:framchangincred0} we are going to investigate the similarity of different frames. 
Using Proposition \ref{sec:perturbframsynth1} and Corollary \ref{sec:upperboundconverg1} the following definition makes sense:
\begin{definition} Let $\mathfrak B_{es} (\Hil)$ be the set of all Bessel sequences in $\Hil$ with the index set $K$. We define the \em Bessel norm \em on this set: For a sequence $(\psi_k) \in \mathfrak B_{es} (\Hil)$ we define the norm
$ \norm{\mathfrak B_{es}}{ ( \psi_k) } = \norm{Op}{D_{(\psi_k)}} $.
\end{definition}
It can be easily shown, that $\norm{\mathfrak B_{es}}{ . }$ is well defined and induces a norm. As shown above it is sufficient to use the operator norm on $c_{00}^2$.
In typical perturbation results like Corollary \ref{sec:perturbBessel1}, for arbitrary sequences it is investigated, if they share the property with another `similar' sequence, which is a Bessel sequence, frame or Riesz basis. In these cases we cannot use the above norm, as we cannot assume at first%the beginning
, that we start out with a Bessel sequence. We are going to define other measures of similarity of sequences, with the property, that if those are small for two Bessel sequences also the Bessel norm is small.

A simple way to measure the similarity of two frames would be in a uniform sense, using %the supremum of 
$\sup \limits_{k} \norm{\Hil}{\psi_k - \phi_k}$, 
but this is not a good measure in general for
frames since it makes an orthonormal basis similar to every
bounded sequence of vectors.
Other similarity measures are more useful, as defined below and motivated in the next result. 
\begin{cor.} \label{sec:framesimi1} Let $(\psi_k)$ be a Bessel sequence, frame, Riesz sequence respectively a Riesz basis. Let $(\phi_k)$ be a sequence with $\sum \limits_k \norm{\Hil}{\psi_k - \phi_k}^2 < A$ (respectively $\sum \limits_k \norm{\Hil}{\psi_k - \phi_k} < A$), then $(\phi_k)$ is a Bessel sequences, frame, Riesz sequence or Riesz basis.

Let $(\phi_k^{(l)} | k \in K)$ be sequences such that for all $\varepsilon$ there exists an $N(\varepsilon)$ with 
$\sum \limits_k \norm{\Hil}{\psi_k - \phi_k^{(l)}}^2 < \varepsilon^2$ (respectively $\sum \limits_k \norm{\Hil}{\psi_k - \phi_k^{(l)}} < \varepsilon$)  for all $l \ge N(\varepsilon)$. Then for $\varepsilon < \sqrt{A}$ and for all $l \ge N_0$ the sequences $(\phi_k^{(l)})$ are Bessel sequences, frames, Riesz sequences respectively Riesz bases with the optimal upper frame bound $B_{opt}^{(l)} \rightarrow B_{opt}$. Furthermore 
$ \norm{Op}{C_{(\phi_k^{(l)})} - C_{(\psi_k)}} < \varepsilon $, 
$ \norm{Op}{D_{(\phi_k^{(l)})} - D_{(\psi_k)}} < \varepsilon $ 
and for $\varepsilon \le 1$ 
$ \norm{Op}{S_{(\phi_k^{(l)})} - S_{(\psi_k})} < \varepsilon \cdot \left( \sqrt{B +1} \cdot \sqrt{B} \right)$.
\end{cor.}
\begin{proof} Let $c \in c_{00}$, then
$$ \norm{\Hil}{D_{(\phi_k)} c - D_{(\psi_k)} c} = \norm{\Hil}{\sum \limits_ k c_k \left( \phi_k - \psi_k \right)} \le \sum \limits_k \left| c_k \right| \norm{\Hil}{\psi_k - \phi_k} \le $$
$$ \sqrt{\sum \limits_k \left| c_k \right|^2} \sqrt{ \sum \limits_k \norm{\Hil}{\psi_k - \phi_k}^2}  \Longrightarrow  \norm{Op}{D_{(\phi_k)} - D_{(\psi_k)}} \le \sqrt{ \sum \limits_k \norm{\Hil}{\psi_k - \phi_k}^2} $$

So in the first case $\norm{Op}{D_{(\phi_k)} - D_{(\psi_k)}} < \sqrt{A}$ and therefore $(\phi_k)$ forms a Bessel sequence, frame, Riesz sequence or Riesz basis. 

In the second case we get $\norm{Op}{D_{\phi_k^{(l)}} - D_{(\psi_k)}} < \varepsilon$ for $l \ge N(\varepsilon)$. With Corollary \ref{sec:upperboundconverg1} the result for the bounds is proved.

$$ \norm{Op}{C_{\phi_k^{(l)}} f - C_{(\psi_k)}} =  \norm{Op}{D_{\phi_k^{(l)}}^* - D_{(\psi_k)}^*} = \norm{Op}{D_{\phi_k^{(l)}} - D_{(\psi_k)}}< \varepsilon $$

$$  \norm{Op}{S_{\phi_k^{(l)}} - S_{\psi_k}} =  \norm{Op}{D_{\phi_k^{(l)}} \circ C_{\phi_k^{(l)}} - D_{(\psi_k)} \circ C_{(\psi_k)} } = $$
$$ =  \norm{Op}{D_{\phi_k^{(l)}} \circ C_{\phi_k^{(l)}} - D_{\phi_k^{(l)}} \circ C_{(\psi_k)} + D_{\phi_k^{(l)}} \circ C_{(\psi_k)} - D_{(\psi_k)} \circ C_{(\psi_k)} } \le $$
$$ \le  \norm{Op}{D_{\phi_k^{(l)}}}  \norm{Op}{C_{\phi_k^{(l)}}  - C_{(\psi_k)}} + \norm{Op}{ D_{\phi_k^{(l)}}  - D_{(\psi_k)}}  \norm{Op}{ C_{(\psi_k)} } $$
$$ \le \sqrt{B +1} \cdot \varepsilon + \varepsilon \cdot \sqrt{B} = \varepsilon \cdot \left( \sqrt{B + 1} + \sqrt{B} \right)  $$
Which follows from Corollary \ref{sec:upperboundconverg1}, as there is an $N(1)$ such that for all $l \ge N(1)$ 
$\norm{Op}{D_{\phi_k^{(l)}}} \le \sqrt{B+1}$. 

For all sequences $C$ we know $\norm{1}{c} \ge \norm{2}{c}$. Therefore the corresponding $l^1$ property above is also true.
\end{proof}

With these similarity measures in general neither a norm nor a metric is defined on the set of Bessel sequences. 
Nevertheless it is useful to use this 'similarity measures' as seen in the last two corollaries. 
\begin{definition} Let $(\psi_k)_{k \in K}$ and $(\psi_k^{(l)})_{k \in K}$ be a sequence of elements for all $l \in \NN$. The sequences $(\psi_k^{(l)})$ are said to \em converge to $(\psi_k)$ in an $l^p$ sense\em , denoted by $(\psi_k^{(l)}) \stackrel{l^p}{\longrightarrow} (\psi_k)$, if $\forall \varepsilon > 0$ there exists an $N > 0$ such that  $\left(\sum \limits_k \norm{\Hil}{\psi_k^{(l)} - \psi_k}^p\right)^{\frac{1}{p}} < \varepsilon$ for all $\ l \ge N$.
\end{definition}
The convergence in an $l^\infty$ sense clearly coincides with uniform convergence, which is not valuable for our purposes, as seen above. The convergence in an $l^1$ and $l^2$ sense is very useful, in contrast, see Section \ref{sec:framchangincred0}.

\section{Bessel Multipliers} \label{sec:frammult0}

In \cite{schatt1} R. Schatten provides a detailed study of ideals of compact operators using their singular value decomposition. He investigates the operators of the form $\sum \lambda_i \varphi_i \otimes_i \overline{\psi_i}$ where $( \varphi_i )$ and $( \psi_i )$ are orthonormal families. We are interested in similar operators where the families are Bessel sequences. 

\begin{definition} 
\label{sec:besselmult1} Let $\Hil_1$, $\Hil_2$ be Hilbert-spaces, let $\left( \psi_k \right)_{k \in K}\subseteq \Hil_1$ and $\left( \phi_k \right)_{k \in K} \subseteq \Hil_2$ be Bessel sequences. Fix $m \in l^\infty (K)$. Define the operator ${\bf M}_{m, ( \phi_k), (\psi_k)} : \Hil_1 \rightarrow \Hil_2$, the \em Bessel multiplier \em for the Bessel sequences $( \psi_k )$ and $( \phi_k )$, as the operator 
$$ {\bf M}_{m, (\phi_k), ( \psi_k )} (f)  =  \sum \limits_k m_k \left< f, \psi_k \right> \phi_k .$$
The sequence $m$ is called the \em symbol \em of $\bf M$. %\index{symbol}
For	frames we will call the resulting Bessel multiplier a \em frame multiplier\em , for Riesz sequence a \em Riesz multiplier\em .
\end{definition}

Let us denote $\MM_{m, (\psi_k)} = \MM_{m, (\psi_k), (\psi_k)}$.
Furthermore let us simplify the notation, if there is no chance of confusion, using $\MM_m$ or even $\MM$. The definition of a multiplier can also be expressed in the following way: 
$$  {\bf M}_{m, (\phi_k) , (\psi_k)} = D_{(\phi_k)} (m \cdot C_{(\psi_k)}) = \sum \limits_k m_k \cdot \phi_k \otimes_i \psi_k $$  

\begin{definition} For fixed Bessel sequences $(\psi_k)$ and $(\phi_k)$, let $\sigma$ be the relation which assigns  the corresponding symbol to a multiplier , $\sigma({\bf M}_{m, (\phi_k), ( \psi_k )}) = m$. %, with $m$ as in Definition \ref{sec:besselmult1}.  
\end{definition}
This relation $\sigma$ does not have to be a well-defined function. This is only the case if the operators $\phi_k \otimes_i \psi_k$ have a basis property, cf. Section \ref{sec:symop1}. 

In \cite{schatt1} the multiplier for orthonormal sequences were investigated and many `nice' properties were shown. A powerful property is the `symbolic calculus' for orthonormal sequences as follows
$$ \MM_{m^{(1)},(e_k)} \circ \MM_{m^{(2)},(e_k)} = \MM_{m^{(1)} \cdot m^{(2)} ,(e_k)} .$$
In the general Bessel sequence case this is not true anymore.

\begin{cor.} \label{sec:combmuleq1}
For two multipliers $\MM_{m^{(1)},(\phi_k),(\psi_k)}$ and $\MM_{m^{(2)},(\zeta_k),(\xi_k)}$ for the Bessel sequences $(\psi_k), (\zeta_k)\subseteq \Hil_2$,$(\phi_k) \subseteq \Hil_2$,$(\xi_k) \subseteq \Hil_1$ and 
$$\MM_{m^{(1)},(\phi_k),(\psi_k)} = \sum \limits_k m^{(1)}_k \left< f, \psi_k \right> \phi_k \mbox{ \  and \  } \MM_{m^{(2)},(\zeta_k),(\xi_k)} =  \sum \limits_l m^{(2)}_l \left< f, \xi_l \right> \zeta_l$$
the combination is % for all $f$
$$ \left( \MM_{m^{(1)},(\phi_k),(\psi_k)} \circ \MM_{m^{(2)},(\zeta_k),(\xi_k)} \right) \left( f \right)= \sum \limits_k \sum \limits_l  m^{(1)}_k m^{(2)}_l \left< f, \xi_l \right> \left<   \zeta_l  , \psi_k \right> \phi_k =
$$
$$ = \left( D_{(\phi_k)} \mathcal M_{m^{(1)}} G_{\psi_k, \zeta_k} \mathcal M_{m^{(2)}} C_{(\xi_k)} \right) (f)  $$
\end{cor.}

Thus in the general Bessel sequence case no exact symbolic calculus can be assumed, i.e. the combination of symbols does not correspond to the combination of the operators. See Section \ref{sec:rieszmult1} for more details on this. 
In general the product of two frame multipliers is not even a frame multiplier any more. 

\subsection{The Multiplier as an Operator from $l^2$ to $l^2$} \label{sec:mulltlt1}

As a preparatory step we will look at this kind of operators on $l^2$. 
Use the symbol $\mathcal M_m$ for the mapping $\mathcal M_m : l^2 \rightarrow l^2$ and $m \in l^p$ (for a $p > 0$) given by the pointwise multiplication $\mathcal M_ m  ( ( c_k ) ) = ( m_k \cdot c_k )$.  
So a Bessel multiplier $\MM_m$ can be written as:
$${\bf M}_m = D \circ \mathcal M_m \circ C$$

As preparation for one of the main results, Theorem \ref{sec:bessmulprop1}, we show:

\begin{lem.} \label{sec:lthsop1} \begin{enumerate}
\item  Let $m \in l^\infty$. The operator $\mathcal M_m : l^2 \rightarrow l^2$ is bounded with $\norm{Op}{\mathcal M_m} = \norm{\infty}{m}$.
\item $\mathcal M_m^* = \mathcal M_{\overline{m}}$
\item Let $m \in l^1$. The operator $\mathcal M_m : l^2 \rightarrow l^2$ is trace class with $\norm{trace}{\mathcal M_m} = \norm{1}{m}$.
\item Let $m \in l^2$. The operator $\mathcal M_m : l^2 \rightarrow l^2$ is a Hilbert-Schmidt ($\HS$) operator with $\norm{\HS}{\mathcal M_m} = \norm{2}{m}$. 
\item Let $m \in c_0$. Then there exist finite sequences $m_N = \left( m_0, \dots, m_N , 0 , 0, \dots \right)$ with $\mathcal M_{m_N} \rightarrow \mathcal M_m$ as operators in  $l^2$. Therefore $\mathcal M_m$ is compact.
\end{enumerate}
\end{lem.}
\begin{proof} 
1.) Clearly
$ \norm{2}{m \cdot c} \le \norm{\infty}{m} \norm{2}{c} \Longrightarrow \norm{Op}{\mathcal M_m} \le \norm{\infty}{m}$.
On the other hand $\mathcal M_m \delta_i = m_i$ $\Longrightarrow \norm{Op}{\mathcal M_m} \ge \norm{\infty}{m}$.

2.) $\left< \mathcal M_m c , d \right>_{l^2} =  \sum \limits_k m_k c_k \cdot \overline{d}_k = \sum \limits_k c_k \cdot \overline{\overline{m}_k d}_k = \left< c , \mathcal M_{\overline{m}} d \right>_{l^2}$.

3.) $\left[ \mathcal M_m \right] = \sqrt{ \mathcal M_m^* \mathcal M_m} = \sqrt{\mathcal M_{\overline{m}}\mathcal M_m} = \mathcal M_{\left| m \right|}$ and so using properties of the trace \cite{wern1} 
$\norm{trace}{\mathcal M_m} = \sum \limits_i  \left< \left[ \mathcal M_m \right] \delta_i , \delta_i \right>  = \norm{1}{m}$.

4.) $\norm{\HS}{\mathcal M_m}^2 = \sum \limits_i \norm{2}{ \mathcal M_m \delta_i} = \norm{2}{m}^2 .$

5.) For $c \in l^2$ %$\lim \limits_{N\rightarrow \infty} m_N \cdot c = m \cdot c $ as 
$\norm{2}{m_N \cdot c - m \cdot c} \le \norm{\infty}{m_N-m} \cdot \norm{2}{c}$ and so 
$\norm{Op}{\mathcal M_{m_N}-\mathcal M_m} \rightarrow 0$.
\end{proof}

\section{Properties of Multipliers} \label{sec:frammulprop}

Equivalent results as proved in \cite{feinow1} for Gabor multiplier can be shown for Bessel multipliers.
\begin{theorem} \label{sec:bessmulprop1} \it Let ${\bf M} = \MM_{m,(\phi_k),(\psi_k)}$ be a Bessel multiplier for the Bessel sequences $(\psi_k) \subseteq \Hil_1$ and $(\phi_k\} \subseteq \Hil_2$ with the bounds $B$ and $B'$. %Let  $\sigma( {\bf M} ) = m$ be its symbol, 
Then
\begin{enumerate}
\item If $m \in l^\infty$  
${\bf M}$ is a well defined bounded operator with 
\mbox{$\norm{Op}{\bf M} \le  \sqrt{B'}  \sqrt{B} \cdot \norm{\infty}{m}$.}
Furthermore the sum $\sum \limits_k m_k \left< f, \psi_k \right> \phi_k$ converges unconditionally for all $f \in \Hil_1$.
\item ${\left( {\bf M}_{m, (\phi_k), (\psi_k)} \right)}^* = {\bf M}_{\overline{m}, (\psi_k), (\phi_k)}$. Therefore if 
$m$ 
is real-valued and $(\phi_k) = (\psi_k)$, ${\bf M}$ is self-adjoint.
\item If $m \in c_0$, 
 $\bf M$ is a compact operator. 
\item If $m \in l^1$, ${\bf M}$ is a trace class operator with $\norm{trace}{M} \le \sqrt{B'} \sqrt{B} \ \norm{1}{m}$. And $\ tr (M) = \sum \limits_k m_k \left< \phi_k , \psi_k \right>$.
\item If $m \in l^2$, ${\bf M}$ is a Hilbert Schmidt operator with $\norm{\mathcal HS}{M} \le \sqrt{B'} \sqrt{B} \ \norm{2}{m}$.
\end{enumerate}
\end{theorem} 

\begin{proof} 
1.) 
$$\norm{Op}{ {\bf M}} = \norm{Op}{C \circ \mathcal M_m \circ D} \le \norm{Op}{C} \cdot \norm{\infty}{m} \cdot \norm{Op}{D} \le \sqrt{B} \norm{\infty}{m} \sqrt{B'} .$$
As $(\phi_k)$ is a Bessel sequence, $\sum c_k \phi_k$ converges unconditionally for all $(c_k) \in l^2$, in particular for $(m_k \cdot \left< f , \psi_k \right>)$. 

2.) ${\bf M}=  C_{(\psi_k)} \circ \mathcal M_m \circ D_{(\phi_k)}  =  C_{(\psi_k)} \circ \mathcal M_m \circ C^*_{(\phi_k)}$, so with Lemma \ref{sec:lthsop1}   ${\bf M}^* = C_{(\phi_k)} \circ M_m^* \circ C_{(\psi_k)}^* = C_{(\phi_k)} \circ M_{\overline{m}} \circ D_{(\psi_k)}$. 

3.) Let $m_N$ be the finite sequences from Lemma \ref{sec:lthsop1}, then 
$$ \norm{Op}{ \MM_{m_N} - \MM_m} = \norm{Op}{D \mathcal M_{m_N} C - D \mathcal M_{m} C} = \norm{Op}{D \left( \mathcal M_{m_N} - \mathcal M_{m} \right) C} \le $$
$$ \le \norm{Op}{D} \norm{Op}{\mathcal M_{m_N} - \mathcal M_{m} } \norm{Op}{C} \le \sqrt{B'} \cdot \varepsilon \sqrt{B} .$$
For every $\varepsilon' = \frac{\varepsilon}{\sqrt{B \cdot B'}}$, there is a $N_\varepsilon$ such that $\norm{Op}{\mathcal M_{m_N} - \mathcal M_{m} } < \varepsilon'$ and therefore $\norm{Op}{ \MM_{m_N} - \MM_m} < \varepsilon$ for all $N \ge N_\varepsilon$.
$\MM_{m_N}$ is a finite sum of rank one operators and so has finite rank. This means that $\MM_m$ is a limit of finite-rank operators and therefore compact.

4.)  
$ {\bf M}(f) = \sum \limits_k \left< f , \psi_k \right> \left( m_k \cdot \phi_k \right)$, 
so according to the definition of trace class operators \cite{wern1} %\ref{sec:trace1} 
we have to show that
$ \norm{trace}{\MM} = \sum_k \norm{\Hil}{\psi_k} \cdot \norm{\Hil}{m_k \phi_k} < \infty$.
$$ \norm{trace}{\MM} = \sum_k \norm{\Hil}{\psi_k} \cdot \norm{\Hil}{m_k \phi_k} = \sum_k \norm{\Hil}{\psi_k} \left| m_k \right| \norm{\Hil}{\phi_k} \le \sqrt{B} \cdot \sqrt{B'} \cdot \norm{1}{m} .$$
$$ \Longrightarrow tr (M) = \sum \limits_k \left< m_k \cdot \phi_k , \psi_k \right> = \sum \limits_k m_k \left< \phi_k , \psi_k \right> .$$

5.) The operator $\mathcal M_m : l^2 \rightarrow l^2$ is in $\HS$ due to Lemma \ref{sec:lthsop1} with bound $\norm{\HS}{\mathcal M_m} = \norm{2}{m}$. Using the properties of $\HS$ operators we get
$ \norm{\HS}{D \mathcal M_m C} \le \norm{Op}{D} \norm{2}{m} \norm{Op}{C} \le \sqrt{B} \sqrt{B'} 
\norm{2}{m}$. 
\end{proof}

For Riesz and orthonormal bases we can show, see Proposition \ref{sec:multrieszopnorm1}%and Theorem \ref{sec:schattth1}
, that if the multiplier is well-defined,
then the symbol must be in $l^\infty$.  This is not true for general Bessel sequences, as can be seen, when using the following frame: 
Let $(e_i | i \in  \NN)$ be an ONB for $\Hil$. Let $\psi_{p,q} = \frac{1}{p} \cdot e_q$. Then $\left( \psi_{p,q} \left| (p,q) \in \NN^2 \right. \right)$
is a tight frame as 
$$ \sum \limits_{p,q} \left| \left< f, \psi_{p,q} \right> \right|^2 = \sum \limits_{p,q} \left| \left< f , \frac{1}{p} \cdot e_q \right> \right|^2 = \sum \limits_{p} \frac{1}{\left|p\right|^2}\sum \limits_{q} \left| \left< f , e_q \right> \right|^2 = \sum \limits_{p} \frac{1}{\left|p\right|^2} \norm{\Hil}{f} = \norm{\Hil}{f} \cdot \frac{\pi^2}{6} $$

Define a symbol $m$ by $m_{p,q} = p^2$. Then 
$$\MM_{m, (\psi_{p,q})} (f) = \sum \limits_{p,q} p^2 \left< f, \frac{1}{p} \cdot e_q \right> \frac{1}{p} \cdot e_q = \sum \limits_{p,q} \left< f, \cdot e_q \right>  \cdot e_q = f .$$
So the operator $\MM_{m_{k,l}, (\psi_{k,l})} = %\equiv
 Id$ is bounded although the symbol is not.

\subsection{From Symbol to Operator} \label{sec:symop1}

When is the operator uniquely defined by the symbol? When is the relation $\sigma$ a function? 
This question is equivalent to the question of whether the sequence of operators $(\phi_h \otimes_i \psi_k)$ forms a Riesz sequence, as the rank one operators $\psi_k \otimes_i \overline{f}_k$ form a Bessel sequence in $\HS$. This follows directly from the following result, as every subsequence of a Bessel sequence is a Bessel sequence again \cite{ole1}. 
\begin{prop.} \label{sec:hsbessel1}Let $(\psi_k, k \in K)$ and $(\phi_k, k \in K)$ be Bessel sequences in $\Hil_1$ respectively $\Hil_2$ with bounds $B_1$ and $B_2$. 
The rank one operators $(\psi_k \otimes_i \overline{\phi}_l)$ with $(k,l) \in K \times K$ form a Bessel sequence in $\HS(\Hil_2,\Hil_1)$ with bounds $B_1 \cdot B_2$.
\end{prop.}
\begin{proof} Let $O \in \HS(\Hil_2,\Hil_1)$. Then by properties of the Hilbert-Schmidt inner product %Corollary \ref{sec:HSinnprodmul1}
$$ \sum \limits_{k,l} \left| \left< O, \psi_k \otimes_i \overline{f}_l \right>_{\HS} \right|^2 = \sum \limits_{l} \sum \limits_{k} \left| \left< O \phi_l, \psi_k \right> \right|^2  \le  B_1 \sum \limits_{l} \norm{\Hil}{O \phi_l}^2 $$
Let now be $(e_i, i \in I)$ be any ONB of $\Hil_2$, then 
$ \sum \limits_{l} \norm{\Hil}{O \phi_l}^2 = \sum \limits_{l} \sum \limits_{i} \left| \left< O \phi_l , e_i \right> \right|^2 = 
\sum \limits_{i} \sum \limits_{l} \left| \left< \phi_l , O^* e_i \right> \right|^2 \le B_2 \cdot \sum \limits_{i} \norm{\Hil_2}{O^* e_i}^2 = $
$ = B_2 \cdot \norm{\HS}{O^*}^2 = B_2 \cdot \norm{\HS}{O}^2 \Longrightarrow \sum \limits_{k,l} \left| \left< O, \psi_k \otimes_i \overline{f}_l \right> \right|^2 \le B_1 B_2 \norm{\HS}{O}^2$.\end{proof}
\section{Riesz Multipliers} \label{sec:rieszmult1}
For Riesz sequences the family $( \psi_k \otimes_i \overline{f}_k)$ is certainly a Riesz sequence in $\HS$, following Proposition \ref{sec:hsbessel1} and the fact that
$$ \left< \psi_k \otimes_i \overline{f}_k, \tilde{g}_l \otimes_i \overline{\tilde{f}}_l \right> = < \psi_k , \tilde{g}_l > \cdot < \tilde{f}_l , \phi_k > = \delta_{k,l} \cdot \delta_{k,l}$$

In this case for $m \in l^2$ the function $m \mapsto \MM_m$ is injective as the multiplier is just the synthesis operator of the sequence $(\psi_k \otimes_i \overline{f}_k)$ applied on $m$ . We can state a more general property:
\begin{lem.} \label{sec:uniqmultRiest1} Let $(\psi_k) \subseteq \Hil_1$ be a Bessel sequence with no zero elements, and $(\phi_k)\subseteq \Hil_2$ a Riesz sequence. Then the mapping $m \mapsto \MM_{m,\phi_k,\psi_k}$ is injective from $l^\infty$ into $\BL(\Hil_1,\Hil_2)$.
\end{lem.}
\begin{proof} 
Suppose $ \MM_{m,(\phi_k),(\psi_k)} =  \MM_{m',(\phi_k),(\psi_k)}$
$ \Longrightarrow \sum \limits_k m_k \left< f , \psi_k \right> \phi_k = \sum \limits_k m'_k \left< f , \psi_k \right> \phi_k$ 
for all $f$. As $\phi_k$ is a Riesz basis for its span
$ \Longrightarrow m_k \left< f , \psi_k \right> = m'_k \left< f , \psi_k \right>$ for all $f, k$. 
For any $k \in K$ we know $\psi_{k} \not=0$. So there exist an $f$ such that $\left< f , \psi_{k} \right> \not= 0$. Therefore
$ m_{k} = m'_{k} \mbox{ for all } k$. And so $(m_k) = (m'_k)$.
\end{proof}

So if the conditions in Lemma \ref{sec:uniqmultRiest1} are fulfilled, the Bessel sequence $( \psi_k \otimes_i \overline{f}_k)$ is a Riesz sequence in $\HS$. 

For Riesz bases the multiplier is bounded if and only if the symbol is bounded:
\begin{prop.} \label{sec:multrieszopnorm1} Let $(\psi_k)$ be a Riesz basis with bounds $A,B$ and ($\phi_k)$ be one with bounds $A',B'$. Then
$$ \sqrt{A A'} \norm{\infty}{m} \le \norm{Op}{\MM_{m,(\phi_k),(\psi_k)}} \le \sqrt{B B'} \norm{\infty}{m} $$
Particularly $\MM_{m,(\phi_k),(\psi_k)}$ is bounded if and only if $m$ is bounded.
\end{prop.}
\begin{proof} Theorem \ref{sec:bessmulprop1} gives us the upper bound.
For the lower bound let $k_0$ be arbitrary, then
\mbox{$\MM_{m,(\phi_k),(\psi_k)} \left( \tilde{\psi}_{k_0} \right) = \sum \limits_k m_k \left< \tilde{\psi}_{k_0} ,\psi_k \right> \phi_k = \sum \limits_k m_k \delta_{k_0,k} \phi_k = m_{k_0} \phi_{k_0} $.} 
Therefore 
$$ \norm{Op}{\MM_{m,(\phi_k),(\psi_k)}} = \sup \limits_{f \in \Hil} \left\{ \frac{\norm{\Hil}{ \MM_{m,(\phi_k),(\psi_k)} \left( f \right) }}{\norm{\Hil}{f}} \right\} \ge \frac{
\norm{\Hil}{\MM_{m,(\phi_k),(\psi_k)} \left( \tilde{\psi}_{k_0} \right)}
}{
\norm{\Hil}{\tilde{\psi}_{k_0}}
} = $$
$ = \frac{\norm{\Hil}{m_{k_0} \phi_{k_0}}}{\norm{\Hil}{\tilde{\psi}_{k_0}}} \ge \frac{\left| m_{k_0} \right| \sqrt{A'}}{\frac{1}{\sqrt{A}}} \ge \sqrt{A' A} \left| m_{k_0} \right|$, 
using Corollary \ref{sec:boundofrieszelem1} and the properties of the dual frame.
\end{proof}
For an orthonormal sequence $(\epsilon_k)$ the combination of multipliers $\MM = \MM_{m,(\epsilon_k)} \cdot \MM_{m',(\epsilon_k)}$ is just the multiplier with symbol $\sigma(\MM) = m \cdot m'$. This is true for all biorthogonal sequences in the following way:
\begin{cor.} \label{sec:commul1} Let $\left( \psi_k \right)$, $\left( \phi_k \right)$, $\left( \zeta_k \right)$ and $\left( \xi_k \right)$ be Bessel sequences, such that $(\phi_k)$ and $(\psi_k)$ are biorthogonal to each other, then 
$$ \left( \MM_{m^{(1)},(\xi_k),(\psi_k)} \circ \MM_{m^{(2)},(\phi_k),(\zeta_k)} \right) \left( f \right)=  \MM_{m^{(1)} \cdot m^{(2)},(\xi_k),(\zeta_k)} $$
\end{cor.}

So for Riesz sequences the we get that symbol of the combination of multipliers is the multiplication of the symbols.
The reverse of this is also true as stated in Corollary \ref{sec:classrieszmul1}. 
For this result we first show the following property:
\begin{prop.} \label{sec:commulequiv1} Let $( \psi_k )$ and $( \phi_k )$ be Bessel sequences in $\Hil$ with the same index set $K$. If $\forall m^{(1)}, m^{(2)} \in c_{00}$ 
$${\bf M}_{m^{(1)}, (\psi_k) , (\phi_k) } \circ {\bf M}_{m^{(2)}, (\psi_k) ,(\phi_k)} = {\bf M}_{m^{(1)} \cdot m^{(2)},(\psi_k) ,(\phi_k)}$$
then for all pairs $(k, l) \in K \times K$ either
$\phi_{l} = 0$, % or
$\psi_k = 0$ %\wedge \tilde{g}_k = 0$ 
or
$\left< \psi_k , \phi_l \right> = \delta_{k,l}$
\end{prop.} 
\begin{proof} Choose $k_0, k_1$ in the index set. Let $m^{(1)} = \delta_{k_0}$ and $m^{(2)} = \delta_{k_1}$. 
$$ {\bf M}_{m^{(1)},(\psi_k) , (\phi_k) } \circ {\bf M}_{m^{(2)}, (\psi_k) , (\phi_k)} = {\bf M}_{m^{(1)} \cdot m^{(2)},(\psi_k) , (\phi_k)}$$ 
is in this case equivalent via Lemma \ref{sec:combmuleq1} to
$$ \left< f , \phi_{k_1} \right> \left< \psi_{k_1} , \phi_{k_0}\right> \cdot \psi_{k_0} = \delta_{k_0 , k_1} \left< f , \phi_{k_1} \right> \psi_{k_0} \  \forall \ f \in \Hil $$

Let $k_1 \not= k_0$ then this means that we obtain
$ \left< f , \phi_{k_1} \right> \left< \psi_{k_1} , \phi_{k_0}\right> \cdot \psi_{k_0} = 0 $.
So either %\begin{enumerate}
$\psi_{k_0} = 0$ or
$\left< f , \phi_{k_1} \right> = 0$ for all $f$ and so $\phi_{k_1} = 0$, or
$\left< \psi_{k_1} , \phi_{k_0}\right> = 0$.

Let $k_1 = k_0$.
$ \left< f , \phi_{k_1} \right> \left(  \left< \psi_{k_1} , \phi_{k_0}\right> - 1 \right) \psi_{k_0} = 0$. 
Either 
$\psi_{k_0} = 0$ or
$\left< f , \phi_{k_1} \right> = 0$ for all $f$ and so $\phi_{k_1} = 0$ or 
$\left< \psi_{k_1} , \phi_{k_0}\right> = 1$. 
\end{proof}

This means we have found a way to classify Riesz bases by multipliers: 
\begin{cor.} \label{sec:classrieszmul1}
Let $( \psi_k )$ and $( \phi_k )$ be Bessel sequences with $\psi_k \not=0$ and $\phi_k \not=0$  for all $k \in K$.  If and only if $\sigma \left( \MM_{m^{(1)},(\phi_k),(\psi_k)} \circ  \MM_{m^{(2)},(\phi_k),(\psi_k)} \right) = \sigma \left(  \MM_{m^{(1)},(\phi_k),(\psi_k)} \right) \cdot \sigma \left(  \MM_{m^{(2)},(\phi_k),(\psi_k)} \right)$ for all multipliers $ \MM_{m^{(1)},(\phi_k),(\psi_k)}$, $ \MM_{m^{(1)},(\phi_k),(\psi_k)}$ with $m^{(1)}$,$m^{(2)}$ finite, then these frames are biorthogonal to each other and therefore have to be  Riesz bases. \rm
\end{cor.}

The commutation of multipliers involving Riesz sequences behaves also very canonically:
\begin{cor.} \label{sec:commriesmul1} Let $\left( \psi_k \right)$ be a Riesz sequence, then 
$$ \MM_{m^{(1)},(\tilde{\psi}_k),(\psi_k)} \circ \MM_{m^{(2)},(\tilde{\psi}_k),(\psi_k)} = \MM_{m^{(2)},(\tilde{\psi}_k),(\psi_k)}  \circ \MM_{m^{(1)},(\tilde{\psi}_k),(\psi_k)} $$
\end{cor.}

Finally we can ask, when a Riesz multiplier is invertible, or more precisely when it is the inverse of another multiplier. Let us call a sequence $(m_k)$ for which $ 0 < \inf \left| m_k \right| \le \sup \left| m_k \right| < \infty $ a \em semi-normalized \em sequence.

\begin{prop.} \label{sec:rieszmultinvert1} Let $(\psi_k)$ and $(\phi_k)$ be Riesz bases and let the symbol $m$ be semi-normalized. 
Then
$\MM_{m_k,(\phi_k),(\psi_k)}^{-1} = \MM_{\frac{1}{m_k},(\tilde{\psi}_k),(\tilde{\phi}_k)} $.
\end{prop.}
\begin{proof} If $(m_k)$ is semi-normalized, $(\frac{1}{m_k})$ is, too. Therefore $(\frac{1}{m_k}) \in l^\infty$. Lemma \ref{sec:combmuleq1} tells us that
$$ \left( \MM_{m,(\phi_k),(\psi_k)} \circ \MM_{\frac{1}{m}, (\tilde{\psi}_k),(\tilde{\phi}_k)} \right) \left( f \right)= \sum \limits_k \sum \limits_l  m_k \frac{1}{m_l} \left< f, \tilde{\phi}_l \right> \left<   \tilde{\psi}_l  , \psi_k \right> \phi_k = $$
$$ = \sum \limits_k \sum \limits_l  m_k \frac{1}{m_l} \left< f, \tilde{\phi}_l \right> \delta_{l,k} \phi_k = \sum \limits_k m_k \frac{1}{m_k} \left< f, \tilde{\phi}_k \right> \phi_k  = f$$
With the commutativity shown in Corollary \ref{sec:commriesmul1} we can finish the proof.
\end{proof} 

\section{Changing the Ingredients} \label{sec:framchangincred0}

A  Bessel multiplier clearly depends on the chosen symbol, analysis and synthesis sequence. A natural question arises: What happens if these items are changed. Are the frame multipliers similar to each other if the symbol or the frames are similar to each other (in the right similarity sense)? 

\begin{theorem} \label{sec:frammulprop1}
Let ${\bf M}$ be a multiplier for the Bessel sequences $(\psi_k)$ and $(\phi_k)$ with Bessel bounds $B_1$ and $B_2$ respectively. Then the operator $\MM$ depends continuously on $m$, $(\psi_k)$ and $(\phi_k)$, in the following sense: Let $(\psi_k^{(l)})$ and $(\phi_k^{(l)})$ be sequences indexed by $l \in\NN$.
\begin{enumerate}
\item 
\begin{enumerate}
\item Let $m^{(l)} \rightarrow m$ in $l^\infty$ then 
$ \norm{Op}{M_{m^{(l)},(\psi_k),(\phi_k)} - M_{m,(\psi_k),(\phi_k)}} \rightarrow 0 $.
\item Let $m^{(l)} \rightarrow m$ in $l^2$ then 
$ \norm{\HS}{M_{m^{(l)},(\psi_k),(\phi_k)} - M_{m,(\psi_k),(\phi_k)}} \rightarrow 0 $.
\item Let $m^{(l)} \rightarrow m$ in $l^1$ then 
$ \norm{trace}{M_{m^{(l)},(\psi_k),(\phi_k)} - M_{m,(\psi_k),(\phi_k)}} \rightarrow 0$.
\end{enumerate}
\item 
\begin{enumerate}
\item Let $m \in l^1$ and let the sequences $(\psi_k^{(l)})$ be Bessel sequences converging uniformly to $(\psi_k)$. 
Then for $l \rightarrow \infty$
$ \norm{trace}{M_{m,(\psi_k^{(l)}),(\phi_k)} - M_{m,(\psi_k),(\phi_k)}} \rightarrow 0$.
\item Let $m \in l^2$ and let the sequences $(\psi_k^{(l)})$ converge to $(\psi_k)$ in an $l^2$ sense. 
Then for $l \rightarrow \infty$
$ \norm{\HS}{M_{m,(\psi_k^{(l)}),(\phi_k)} - M_{m,(\psi_k),(\phi_k)}} \rightarrow 0$.
\item Let $m \in l^\infty$ and let the sequences $(\psi_k^{(l)})$ converge to $(\psi_k)$ in an $l^1$ sense. 
Then for $l \rightarrow \infty$
$ \norm{Op}{M_{m,(\psi_k^{(l)}),(\phi_k)} - M_{m,(\psi_k),(\phi_k)}} \rightarrow 0$.
\end{enumerate}
\item For Bessel sequences $(\phi_k^{(l)})$ converging to $(\phi_k)$,
corresponding properties as in (2) apply.
\item 
\begin{enumerate}
\item Let $m^{(l)} \rightarrow m$ in $l^1$, $(\psi_k^{(l)})$ and $(\phi_k^{(l)})$ be Bessel sequences with bounds  $B_1^{(l)}$ and $B_2^{(l)}$, such that there exists $\bf B_1$ and $\bf B_2$ with $B_1^{(l)} \le \bf B_1$  and $B_2^{(l)} \le \bf B_2$. Let the sequences $(\psi_k^{(l)})$ and $(\phi_k^{(l)})$ converge uniformly to $(\psi_k)$ respectively $(\phi_k)$.
Then for $l \rightarrow \infty$
$ \norm{trace}{M_{m^{(l)},(\psi_k^{(l)}),(\phi_k^{(l)})} - M_{m,(\psi_k),(\phi_k)}} \rightarrow 0$.
\item Let $m^{(l)} \rightarrow m$ in $l^2$ and let the sequences $(\psi_k^{(l)})$ respectively $(\phi_k^{(l)})$ converge to $(\psi_k)$ respectively $(\phi_k)$ in an $l^2$ sense. 
Then for $l \rightarrow \infty$
$ \norm{\HS}{M_{m,(\psi_k^{(l)}),(\phi_k)} -M_{m,(\psi_k),(\phi_k)}} \rightarrow 0$.
\item Let $m^{(l)} \rightarrow m$ in $l^\infty$ and let the sequences $(\psi_k^{(l)})$ respectively $(\phi_k^{(l)})$ converge to $(\psi_k)$ respectively $(\phi_k)$ in an $l^1$ sense.
Then for $l \rightarrow \infty$
$ \norm{Op}{M_{m,(\psi_k^{(l)}),(\phi_k)} - M_{m,\psi_k,\phi_k}} \rightarrow 0$.
\end{enumerate}
\end{enumerate}
\end{theorem} 
\begin{proof} 
1.) For a sequence of symbols this is a direct result of Theorem \ref{sec:bessmulprop1} and
$\norm{\HS}{\MM_{m^{(l)},(\psi_k),(\phi_k)} - \MM_{m,(\psi_k),(\phi_k)}} = \norm{\HS}{\MM_{\left( m^{(l)} - m\right),(\psi_k),(\phi_k)}} \le $
$ \norm{2}{m^{(l)}-m} \sqrt{B B'}$. 
The result for the operator and infinity norm respectively trace and $l^1$ norms can be proved in an analogue way. 

2.) For points (b) and (c) we know from Corollary \ref{sec:framesimi1} that the sequences are Bessel sequences. 
For all the norms $(Op,\HS,trace)$ $\norm{}{\psi_k \otimes_i \phi_k} = \norm{\Hil}{\psi_k} \norm{\Hil}{\phi_k} $ and so
$ \norm{}{\sum m_k \psi_k^{(l)} \otimes_i \phi_k - \sum m_k \psi_k \otimes_i \phi_k} =  \norm{}{\sum m_k \left( \psi_k^{(l)}- \psi_k \right) \otimes_i \phi_k} \le $
$ \sum \limits_k \left| m_k \right| \norm{\Hil}{\psi_k^{(l)}- \psi_k} \sqrt{B'} = (*)$ \\
$ \mbox{ case a : } (*) \le \sqrt{B'} \left( \sum \limits_k \left| m_k \right| \right) \sup \limits_l \left\{ \norm{\Hil}{\psi_k^{(l)}- \psi_k} \right\} \le \sqrt{B'} \norm{1}{m} \varepsilon  $ \\
$ \mbox{ case b : } (*) \le \sqrt{B'} \sqrt{ \sum \limits_k \left| m_k \right|^2 } \sqrt{ \sum \norm{\Hil}{\psi_k^{(l)}- \psi_k}^2 } \le \sqrt{B'} \norm{2}{m} \varepsilon $ \\
$ \mbox{ case c : } (*) \le \sqrt{B'} \norm{\infty}{m} \sum \norm{\Hil}{\psi_k^{(l)}- \psi_k} \le \sqrt{B'} \norm{\infty}{m} \varepsilon $

3.) Use a corresponding  argumentation as in (2).

4.) For points (b) and (c) Corollary \ref{sec:upperboundconverg1} states that $(\psi_k^{(l)})$ and $(\phi_k^{(l)})$ are Bessel sequences and there are common Bessel bounds $\bf B_1$ and $\bf B_2$ for $l \ge N_1$. So using the results above we get
$$ \norm{}{M_{m^{(l)},(\psi_k^{(l)}),(\phi_k^{(l)})} - M_{m,(\psi_k),(\phi_k)}} \le  \norm{}{M_{m^{(l)},(\psi_k^{(l)}),(\phi_k^{(l)})} - M_{m,(\psi_k^{(l)}),(\phi_k^{(l)})}} + $$
$$ + \norm{}{M_{m,(\psi_k^{(l)}),(\phi_k^{(l)})} - M_{m,(\psi_k),(\phi_k^{(l)})}} + \norm{}{M_{m,(\psi_k),(\phi_k^{(l)})} - M_{m,(\psi_k),(\phi_k)}} \le $$
$$ \le \varepsilon \sqrt{\bf B B'} + \norm{}{m} \varepsilon \sqrt{\bf B'} + \norm{}{m} \sqrt{B} \varepsilon = \varepsilon \cdot \left( \sqrt{\bf B B'} + \norm{}{m} \left( \sqrt{\bf B'} + \sqrt{B} \right) \right) $$
for an $l$ bigger than the maximum $N$ needed for the convergence conditions. This is true for all pairs or norms $(Op, \infty)$, $(\HS,l^2)$ and $(trace,l^1)$.
\end{proof}

\section{Perspectives} \label{sec:conclpersp0}

For the future many questions are still open. For example it seems very likely, that for symbols $m \in l^p$ the multiplier lies in the $p$-Schatten operator class. Connected to that an investigation of the singular values of these operators might be worthwhile. 
The combination of two multipliers are connected to the Gram matrix, see Lemma \ref{sec:combmuleq1}. It will be interesting to apply the results for the decay properties of the Gram matrix in \cite{forngroech1} to this topic. An interesting step away from the unstructured frame will be to investigate frame multipliers for structured frames, found e.g. in \cite{yonboek1}.
The topic of frame multipliers is closely related to the notion of \em weighted frames \em as introduced in \cite{bogdvan1}. 
It can be easily proved that for a positive, semi-normalized symbol the multiplier corresponds to the frame operator of a weighted frame. 
This connection should be investigated further and the theory of frame multiplier should be applied to the further context of the paper \cite{bogdvan1}, computational issues for wavelets on the sphere.

Applications of these objects already exist. It seems that acoustics is a very interesting field for that. Frame multipliers there are not only used as (irregular or regular) Gabor multipliers like in \cite{xxlphd1} or \cite{majxxl1}, but also as a multipliers for a gammatone filter bank in \cite{piche1}. In the engineering literature frame multipliers for regular Gabor frames are known as \em Gabor filters \em \cite{hlawatgabfilt1}. Also first ideas are investigated to use this concept with wavelets to apply in the context of evaluations of 
noise barriers. The importance of the theoretical results in these and other application should be investigated further.

\section*{Acknowledgments}

The author would like to thank Hans G. Feichtinger, Bruno Torr\'esani, Jean-Pierre Antoine and Wolfgang Kreuzer for
many helpful comments and suggestions. He would like to thank the hospitality of the LATP, CMI, Marseille, France and FYMA, UCL, Louvain-la-Neuve, Belgium, where part of this work was prepared, supported by the HASSIP-network.

\end{document}